\documentclass[12pt]{article}

\usepackage{a4wide}
\usepackage{amssymb}
\usepackage{amsfonts}
\usepackage{amsmath}

\date{}

\newtheorem{proposition}{Proposition}[section]
\newtheorem{theorem}[proposition]{Theorem}
\newtheorem{lemma}[proposition]{Lemma}

\newtheorem{corollary}[proposition]{Corollary}
\newcommand{\smallbox}{{\vrule height3pt width3pt depth0pt}}

\def\GK{{\rm  GK}\,}
\def\Kdim{{\rm K.dim }\,}

\def\der{\partial }

\def\nFM0{{\nu }_{F,M_0}}
\def\nFN0{{\nu }_{F,N_0}}
\def\nGN0{{\nu }_{G,N_0}}

\def\N0{ {\bf N}_0 }

\def\t{\otimes}
\def\g{\gamma}
\def\v{\varphi}
\def\ra{\rightarrow}
\def\lra{\leftrightarrow}
\def\Xpm{X^{\pm }}

\def\s{\sigma}

\def\l1{{\lambda}_1}

\def\a{\alpha}
\def\a0{ {\alpha }_0}
\def\a1{ {\alpha }_1}

\def\l{\lambda}

\def\nFGM0{{\nu }_{F,G,M_0}}

\def\nFN0{{\nu}_{F,N_0}}

\def\sn{{\sigma}^n}
\def\sm{{\sigma}^m}
\def\si{{\sigma}^i}

\def\sm1{{\sigma}^{-1}}

\def\smtp1{{\sigma}^{-t+1}}

\def\S1{S^{-1}}

\def\Xpm1{X^{\pm 1}_1}

\def\sPM1{{\sigma }^{\pm 1}}
\def\sMP1{{\sigma }^{\mp 1 }}

\def\d{\delta}

\def\CA{{\cal A}}

\def\Ytm1{Y^{t-1}}
\def\Yim1{Y^{i-1}}

\def\CG{{\cal G}}

\def\Aut{{\rm Aut}}
\def\bK{\bar{K}}

\def\bA{\bar{A}}

\def\ad{{\rm ad }}
\def\dim{{\rm dim }}

\def\ker{ {\rm ker } }
\def\ndeg{ {\rm ndeg} }
\def\CJ{ {\cal J}}
\def\gr{ {\rm gr} }
\def\gcd{ {\rm gcd } }
\def\D{ \Delta }
\def\Ev{ {\rm Ev } }

\begin{document}

\author{V. V.  Bavula }

\title{Dixmier's Problem 5 for the Weyl Algebra}

\maketitle
\begin{abstract}
In the paper \cite{Dix}, J. Dixmier posed six problems for the
first Weyl algebra.
In this paper we give a solution to the
Dixmier's Problem 5 from this paper. Problem 3 was solved by
Joseph and Stein \cite{josclA1} (using results of McConnel and
Robson \cite{MR}). Using a (difficult) polarization theorem for
the first Weyl algebra  Joseph \cite{josclA1} solved problem 6 (a
short proof to  this problem is given in \cite{BavDP6}, note that
the same result and the proof are true for the ring of
differential operators on an arbitrary smooth irreducible
algebraic curve \cite{BavDP6}). Problems 1, 2, and 4 are still
open.

\end{abstract}


\section{Introduction}
 Let $K$ be a field of characteristic
zero. The {\em first Weyl algebra} $A_1$
 is an associative $K$-algebra generated over $K$ by elements $X$ and $Y$ subject to
the defining relation $YX-XY=1$. The $n$'th {\em Weyl algebra} $A_n$ is the tensor
 product $A_1\t \cdots \t A_1$ of $n$ copies of the first Weyl algebra. The Weyl
algebra $A_n$ is a simple Noetherian domain of Gelfand-Kirillov dimension $2n$ which
is canonically isomorphic to the ring of differential operators
 $K[X_1, \ldots , X_n, \frac{\der}{\der X_1}, \ldots , \frac{\der}{\der X_n}]$
 with polynomial coefficients. The Weyl algebras have been intensively studied
during the last fifty years. The Gelfand-Kirillov dimension and the transcendence
 dimension of the Weyl algebra $A_n$ were computed by Gelfand and Kirillov,
\cite{GK}. The fact that each derivation of the Weyl algebra $A_n$ is an inner
derivation
 was proved by Dixmier, \cite{Dix3}. The commutativity of the centralizer of an
 arbitrary nonzero
element of the first Weyl algebra was proved by Amitsur, \cite{Ami}. The structure
of maximal commutative subalgebras of the first Weyl algebra was studied and the generators
 of the group of all algebra automorphisms of $A_1$ were found by Dixmier, \cite{Dix}
 (see \cite{Bav-Jor} for a generalization of these results to noncommutative
deformations of type-$A$ Kleinian singularities).
 Simple $A_1$-modules were classified by Block, \cite{Bl}
(see also \cite{Bav} for an alternative approach and some generalizations). The global dimension of $A_n$
 is $n$, this was proved by Rinehart, \cite{Ri}  in the case $n=1$, and by Roos
in the general case. Rentschler and Gabriel proved that the Krull
dimension of $A_n$ is $n$, \cite{RG}. The finite dimensionality of
the vector spaces  ${\rm Ext}^i_{A_1}$ and ${\rm Tor}^{A_1}_i$ for
simple $A_1$-modules was established by McConnell and Robson,
\cite{MR}. The fact that the Gelfand-Kirillov dimension of a
nonzero finitely
 generated $A_n$-module is not less than $n$ ({\it the Bernstein Inequality}) was
proved by Bernstein, \cite{Be}. A finitely generated $A_n$-module of Gelfand-Kirillov dimension $n$
is called a {\em holonomic} module. Each simple module over the first Weyl algebra is
holonomic. The situation is completely different for the Weyl algebras $A_n$, $n\geq 2$.
 The first examples of {\em non-holonomic} $A_n$-modules were constructed by Stafford,
\cite{St},  futher progress in this direction was made by
Coutinho, \cite{Cou}.
 Bernstein and Lunts, \cite{Be-Lu} in the case of the second Weyl algebra, and Lunts
 in the general case, \cite{Lu}, showed that ``generically'' a simple $A_n$-module is
non-holonomic and has Gelfand-Kirillov dimension $2n-1$.  Simple holonomic $A_2$-modules
 were classified by the author and van Oystaeyen, \cite{BVO}. Skew subfields of the
$n$'th Weyl skew field which are invariant under the action of a finite
group were studied by Alev and Dumas, \cite{Ale-Dum}. Makar-Limanov,
 \cite{Mak-Lim},  proved that the first Weyl skew field contains a free
  subalgebra.

In his fundamental paper \cite{Dix} Dixmier initiated a systematic
study of the structure of the first Weyl algebra $A_1$. At the end
of his paper he
 posed 6 problems. In this paper a {\em negative} answer is given to the 5'th problem. {\em Problem 1}
 concerns the question {\em whether an algebra endomorphism of $A_1$
is an algebra automorphism?} A positive answer to a similar problem but for the $n$'th
 Weyl algebra implies the {\em Jacobian Conjecture} as was shown by  Bass,
 Connel and  Wright, \cite{BCW}. For an arbitrary non-scalar element $u$
 of  $A_1$ one can associate the {\em inner derivation} $\ad \, u$ of the
 Weyl    algebra $A_1$, $\ad \,u (a)=ua-au$, $a\in A_1$, and then the $\mathbb{N}$-filtered
 algebra $N(u)=\cup_{i\geq 0}\, N(u, i)$ where $N(u, i):= \ker  (\ad \, u)^{i+1}$.
The zero component of this filtration, $\ker \, \ad \, u$, is the {\em centralizer}
 $C(u)$ of the element $u$ in $A_1$. The algebra $C(u)$ is a commutative algebra which is
a free finitely generated module over its polynomial subalgebra
$K[u]$, \cite{Ami} and \cite{Dix}. Dixmier partitioned all
non-scalar elements of the Weyl algebra $A_1$ into 5 classes $\D_1
, \ldots , \D_5$, and classified
 up to the action of the group $\Aut_K(A_1)$ elements from the class $\D_1$,
 so-called {\em elements of strongly nilpotent type}, and elements from the class
 $\D_3$, so-called {\em elements of strongly semi-simple type}. Problems 2-6 are concerned with
 properties and classification of elements from the remaining classes $\D_2$, $\D_4$
 and $\D_5$. A non-scalar element $u\in A_1$ with $C(u)\neq N(u)$ (resp.,
with $N(u)=A_1$) is called an element of {\em nilpotent type}
(resp., {\em  of strongly nilpotent type}). A non-scalar element
$u\in A_1$ of nilpotent type belongs to $\D_2$ iff $C(u)\neq N(u)\neq A_1$,
and we say that $u$ is of {\em weakly nilpotent type}.

{\bf Dixmier's Problem 5, \cite{Dix}}: {\em Let $u\in A_1$ be an
element of nilpotent type. Set
 $I_n=(\ad \, u)^n N(u,n)$; this an ideal of $C(u)$. Is $I_{n+1}=I_1I_n$ for
$n$ sufficiently large?}

The next Theorem shows that the answer in general is negative. This result is proved
in Section 4, Corollary \ref{Dixpr}

For given natural numbers $n$ and $m\neq 0$, there exist and
unique natural numbers $l$ and $r$ such that $n=lm+r$ and $0\leq
r<m$. The number $l$  is denoted by $[\frac{n}{m}]$.

\begin{theorem}\label{Dixprint} 
Let $\alpha (H)\in K[H]$ be a polynomial of degree $d\geq 1$ in the variable $H=YX$.
The centralizer of the element $u=\alpha (YX)X\in A_1$ is the polynomial ring
 $K[u]$, and $I_k=u^{k-[\frac{k}{d+1}] }\, K[u]$,  for all $k\geq 1$. In particular,
$I_1=uK[u]$ and $I_{i(d+1)-1}=I_{i(d+1)}=u^{id}K[u]$, for all
$i\geq 1$. Hence, $I_1I_{i(d+1)-1}\neq I_{i(d+1)}$, for all $i\geq
1$; and so the Dixmier's Problem 5
 has negative answer.
\end{theorem}

In order to prove this result we consider the Weyl algebra as the
{\em generalized Weyl algebra} $A_1=K[H](\s , H)=\oplus_{i\in
\mathbb{Z}} \, K[H]v_i$ (see Section 2 for details). The
localization of $A_1$ at the Ore set $S=K[H]\backslash \{ 0\}$ is
the skew Laurent extension
 $B=K(H)[X,X^{-1}; \s ]=\oplus_{i\in \mathbb{Z}}\, K(H)X^i $ with the $K$-automorphism $\s $
of the field of rational functions $K(H)$ defined by  $\s
(H)=H-1$. The Weyl algebra $A_1$ is a homogeneous subalgebra of
the $\mathbb{Z}$-graded algebra $B$. In Section 2,
 for a homogeneous element $\alpha X^i$, $\alpha \in K(H)$, $i\in \mathbb{Z}$, the centralizer
 $C(u,B)$  (Proposition \ref{CUB})  and  the algebra $N(u,B)$ (Theorem \ref{Nhomel})
are described.  In Section 3, using these results,
for an arbitrary homogeneous element $u$ of the Weyl algebra $A_1$,  the centralizer
 $C(u,A_1)$  (Proposition \ref{centralizerhomweyl})  and  the algebra $N(u,A_1)$ (Theorem \ref{Nhomweyl})
 are found. Then, for the element
 $u=\alpha X$ as in the Theorem \ref{Dixprint}, we can describe the algebra $N(u, A_1)$ and the
ideals $I_n$. In Section 5, we classify homogeneous elements of
the Weyl algebra $A_1$
 with respect to the Dixmier partition of elements of $A_1$ into the classes $\D_i$.
We prove (Corollary \ref{Dixpr4hom}) that for an arbitrary
homogeneous element $u$ of the Weyl algebra $A_1$ of nilpotent
type : (i) {\em the Dixmier's Problem 4} has positive answer, that
is the associated graded algebra $\CG (u, A_1)$ of the
$\mathbb{N}$-algebra  $N(u,A_1)$ is an affine commutative algebra,
hence Noetherian, and as a consequence the algebra $N(u,A_1)$ is
affine Noetherian; (ii) the Weyl algebra $A_1$ is {\em not} a
finitely generated (left and right) $N(u, A_1)$-module.

For more information about the Weyl algebras the reader is referred
 to the following books \cite{Bj, KL, MRb}.


\section{Centralizer of a Homogeneous Element of the Algebra $B$}

Let $D$ be a ring with an automorphism $\s $ and a central element $a$.
 The {\bf generalized Weyl algebra} $A=D(\sigma, a)$  of degree 1,
is the ring generated by $D$ and two indeterminates $X$ an $Y$  subject to
the relations:
$$
X\alpha=\sigma(\alpha)X \ {\rm and}\ Y\alpha=\sigma^{-1}(\alpha)Y,\; {\rm for \; all }\;
\alpha \in D, \ YX=a \ {\rm and}\ XY=\sigma(a).
 $$
The algebra $A={\oplus}_{n\in \mathbb{Z}}\, A_n$ is a
$\mathbb{Z}$-graded algebra where $A_n=Dv_n$,
 $v_n=X^n\,\, (n>0), \,\,v_n=Y^{-n}\,\, (n<0), \,\,v_0=1.$
 It follows from the above relations that
$$v_nv_m=(n,m)v_{n+m}=v_{n+m}<n,m> $$
for some $(n,m)=\s^{-n-m}(<n,m>)\in D$. If $n>0$ and $m>0$ then
$$n\geq m:\, (n,-m)=\sigma^n(a)\cdots \sigma^{n-m+1}(a),\, (-n,m)=\sigma^{-n+1}(a)\cdots \sigma^{-n+m}(a),$$
$$n\leq m:\,\,\,\, (n,-m)=\sigma^{n}(a)\cdots \sigma(a),\,\,\,
(-n,m)=\sigma^{-n+1}(a)\cdots a,$$
in other cases $(n,m)=1$.

Let $K[H]$ be a polynomial ring in one variable $H$ over the field $K$, $\s :H\ra H-1$ be the
$K$-automorphism of the algebra $K[H]$ and $a=H$. The first Weyl algebra
$A_1=K<X, Y\;|\;YX-XY=1>$ is isomorphic to the generalized Weyl algebra
$$A_1\simeq K[H](\s , H),\; X\lra X,\; Y\lra Y,\; YX\lra H.$$
 We identify both these algebras via this isomorphism, that is
$A_1=K[H](\s , H)$ and $H=YX$.

If $n>0$ and $m>0$ then
$$n\geq m:\; (n,-m)= (H-n)\cdots (H-n+m-1),\; (-n,m)=(H+n-1)\cdots (H+n-m),$$
$$n\leq m:\; (n,-m)=(H-n)\cdots (H-1), \; (-n,m)=(H+n-1)\cdots H,$$
in other cases $(n,m)=1$.

The localization $B=S^{-1}A_1$ of the Weyl algebra $A_1$ at the Ore subset
 $S=K[H]\backslash \{ 0\}$ of $A_1$ is
the {\em skew Laurent polynomial ring} $B=K(H)[X, X^{-1}; \s ]$ with
coefficients from the field $K(H)=S^{-1}K[H]$ of rational functions
and $\s \in \Aut_K \, K(H)$,  $\s (H)=H-1$. The map $A_1\ra B$, $a\ra
a/1$,
is an algebra monomorphism. We identify the algebra $A_1$ with its image in the algebra
$B$, in more detail, via the algebra monomorphism
$$A_1\ra B, \;\; X\ra X, \;\; Y\ra HX^{-1}.$$
 The subalgebra $\CA :=K[H][X,X^{-1}; \s ]$ of $B$ contains the Weyl
 algebra $A_1$ (since the algebra generators $X$ and $Y=HX^{-1}$ of $A_1$ belong to
$\CA $), moreover, $\CA $ is the localization $\CA =A_{1,X}$ of $A_1$ at the
powers of the element $X$. Clearly, $B=S^{-1}\CA $.
 The algebras $B=\oplus_{i\in \mathbb{Z}}\, B_i$ and $\CA =\oplus_{i\in
\mathbb{Z}}\, \CA_i$ are $\mathbb{Z}$-graded algebras where
$B_i=K(H)X^i$ and
 $\CA_i=K[H]X^i$. The algebras $A_1$ and $\CA $ are $\mathbb{Z}$-graded subalgebras
of $B$.

A polynomial $f(H)=\l_nH^n+\l_{n-1}H^{n-1}+\cdots +\l_0\in K[H]$
of degree $n$ is called a {\em  monic}  polynomial  if the leading
coefficient  $\l_n$ of $f(H)$ is $1$. A rational function $h\in
K(H)$ is called a {\em monic}  rational function  if $h=f/g$ for
some monic polynomials $f,g$. A homogeneous element $u=\alpha x^n$
of $B$ is called {\em monic} iff
 $\alpha $ is a monic rational function. We can extend in the obvious way the notion
of degree of a polynomial to the field of rational functions setting,
 $\deg_H \, h = \deg_H\, f - \deg_H\, g$, for $h=f/g\in K[H]$. If $h_1, h_2 \in K(H)$
then $\deg_H \, h_1h_2=\deg_H\,h_1 +\deg_H\, h_2$, and
 $\deg_H (h_1+h_2)\leq \max \{ \deg_H\,h_1 , \deg_H\, h_2\}$.
  We denote by ${\rm sign} (n)$ and by $|n|$ the sign and the
absolute value of $n\in \mathbb{Z}$, respectively.

\begin{proposition}  
\label{CUB} {\sc (Centralizer of a Homogeneous Element of the Algebra $B$)}
\begin{enumerate}
\item Let $u=\alpha X^n$ be a monic element of $B_n$  with $n\neq 0$.  The centralizer $C(u,B)=K[v,v^{-1}]$ is a Laurent polynomial
ring in a uniquely defined variable $v=\beta X^{{\rm sign} (n) s}$
where $s$ is the minimal positive divisor of $n$ for which there
exists an element $\beta =\beta_s \in K(H)$, necessarily monic and
uniquely defined, such that
\begin{equation}
\label{bap}
\beta \,\s^s (\beta )\, \s^{2s}(\beta ) \cdots \s^{(n/s -1)s}(\beta )=\alpha, \;\; {\rm if}\;\; n>0,
\end{equation}
\begin{equation}
\label{bam}
\beta \, \s^{-s} (\beta ) \,\s^{-2s}(\beta ) \cdots \s^{-(|n|/s -1)s}(\beta )=\alpha, \;\; {\rm if}\;\; n<0.
\end{equation}

\item Let $u\in K(H)\backslash K$. Then $C(u,B)=K(H)$.
\end{enumerate}
\end{proposition}

{\it Proof}. The element  $u$ is a homogeneous element of the
$\mathbb{Z}$-graded algebra  $B$, hence its centralizer
$$C=C(u,B)=\oplus_{i\in \mathbb{Z}}\, C_i, \;\; C_i=C\cap B_i,$$
is a graded subalgebra of $B$. Consider
$$H\equiv H(u,B):=\{ i\in \mathbb{Z}\,|\, C_i\neq 0 \}.$$
The set $H$ is a subgroup of $\mathbb{Z}$: $0\in H$, since $1\in
C$;
 $H+H\subseteq H$, since $C_iC_j\subseteq C_{i+j}$ and $B$ is a
domain; $-H\subseteq H$, since if $0\neq v\in C_i$, then $v^{-1}\in
C_{-i}$.

1. In this case, $H=\mathbb{Z}s$ for a uniquely defined positive
divisor $s$ of $n$ $(n\in H$, since $u\in C_n)$.

{\it Claim} 1: {\em for every} $i\in H$, $C_i=K^*\alpha_iX^i$ {\em
for a uniquely defined monic } $\alpha_i\in K(H)$. Obviously,
$\beta X^i\in C_i$, for some $0\neq \beta \in K(H)$, iff
\begin{equation}  
\label{sisn} \frac{\si (\alpha )}{\alpha } =\frac{\sn (\beta )}{
\beta }.
\end{equation}
So, if $0\neq \beta_1X^i,   \beta_2X^i\in C_i$, then
$$ \frac{\sn (\beta_1 )}{\beta_1}= \frac{\sn (\beta_2 )}{ \beta_2 }$$
or, equivalently, $\sn (\beta_2/\beta_1)=\beta_2/\beta_1\in K(H)^{\sn
}=K^*$, this finishes the proof of the claim.

It follows from the claim and from $H=\mathbb{Z}s$, that
$C(u,B)=K[v,v^{-1}]$ for some $v$.

By the claim 1, there exists a unique monic element $0\neq \beta \in
K(H)$ such that $v=\beta X^{{\rm sign} (n) s}$. If $n>0$, then
$$C_n=K^*\alpha X^n\ni v^{n/s}=(\beta X^s)^{n/s}=
\beta \,\s^s (\beta ) \cdots \s^{(n/s -1)s}(\beta )X^n\equiv
\beta_nX^n,$$
hence $\alpha =\beta_n$,
since $\beta_n$ is a monic polynomial.
 If $n<0$, then
$$C_n=K^*\alpha X^n\ni (\beta X^{-s})^{|n|/s}=
\beta \, \s^{-s} (\beta ) \,\s^{-2s}(\beta ) \cdots \s^{-(|n|/s -1)s}(\beta )
X^n\equiv \beta_nX^n,$$
hence $\alpha =\beta_n$,
since $\beta_n$ is a monic polynomial.

{\it Claim} 2: {\em suppose that for some positive divisor} $s$ {\em
of } $n$  {\em and for some } $0\neq \beta \in K(H)$ {\em one of the
corresponding equalities}, (\ref{bap}) {\em or} (\ref{bam}), {\em holds. Then }
$\beta  X^{{\rm sign} (n) s}\in C$. Consider the case $n>0$. Then

$$ \frac{\s^s (\alpha )}{ \alpha }=\s^s (\beta )\, \s^{2s}(\beta ) \cdots
\s^{(n/s)s}(\beta )/\beta \,\s^s (\beta )\, \s^{2s}(\beta ) \cdots
\s^{(n/s -1)s}(\beta )=\frac{ \s^n (\beta )}{ \beta},$$ hence
$\beta X^s\in C$, by (\ref{sisn}). Claim 2 proves the minimality
of the $s$ (in the Proposition).

2. The  centralizer $C$ of $u$ is a homogeneous subalgebra of $B$ which contains $K(H)$.
 A homogeneous element $\beta X^i$ of $B$ with $i\neq 0$ commutes with $u$ iff $\beta =0$ since
$0=[\beta X^i, u]=\beta (\si (u)-u)X^i$ and $u\not\in K=K(H)^{\si }$. This proves that
$C=K(H)$. {\smallbox}

{\sc Definition}. The uniquely defined element $v$ from Proposition \ref{CUB}.(1) is
called the {\em canonical generator} of the algebra $C(u,B)$.

The set of polynomials
\begin{equation}  
\label{polynomialsvn} \v_0 :=1, \; \v_n:=(-1)^n \frac{
H(H+1)\cdots (H+n-1)}{ n! }=(-1)^n \frac{ H\s^{-1} (H)\cdots
\s^{-n+1} (H)}{ n! }, \; n\geq 1,
\end{equation}
is a $K$-basis of the polynomial algebra $K[H]$, $\deg \, \v_n=n$ and
$$ \s (\v_n)-\v_n=\v_{n-1}, \; {\rm for \; all } \; n\geq 0, \; \v_{-1}:=0.$$


Let $K$ be a field and let $K[t]$ be a polynomial ring in an indeterminate
 $t$. Let $M$ be a $K[t]$-module. For an element $p\in K[t]$ we denote by
 $\ker \, p_M$ the kernel of the $K$-linear map $p=p_M:M\ra M$, $m\ra pm$.
 The kernel $\ker \, p_M$ is a $K[t]$-submodule of $M$. For $i\geq 0$, let
$N_i =N_i(t, M):=\ker \, t^{i+1}_M$. Then
$$ N_0\subseteq N_1\subseteq \cdots \subseteq N_i\subseteq \cdots , \;
tN_0=0, \;\; {\rm and } \;\; tN_j\subseteq N_{j-1}, \;\; {\rm for } \;\; j\geq 1.$$
Clearly, $N=N(t,M):=\cup_{i\geq 0}\, N_i$ is a $K[t]$-submodule of $M$.
 We set $N_{-1}=0$. If $0\neq u\in N$ then the unique $i$ such that
 $u\in N_i\backslash N_{i-1}$ is called the {\em nilpotent degree} of $u$,
denoted by $\ndeg \, u$.

{\sc Definition}. A $K[t]$-module $M$ is called a {\em Jordan } $K[t]$-module
iff $M=N(t,M)$ and $tM=M$.

If $M$ is a nonzero Jordan $K[t]$-module then $N_i\neq N_{i+1}$ for
all $i\geq 0$ since otherwise $M=N_j$ for some $j$, hence $M=t^{j+1}M=0$, a contradiction.
 From this fact we conclude that each nonzero Jordan $K[t]$-module is not
 a finitely generated module (hence, is not noetherian).

{\it Example}. A vector space $\CJ =\oplus_{i\geq 0} \, Ke_i$ with
the $K[t]$-module structure defined by $te_0=0$ and
$te_i=e_{i-1}$, $i\geq 1$,  is a Jordan module. The module $\CJ $
is isomorphic to the $K[t]$-module $K[t, t^{-1}]/K[t]$, and
$\CJ_{i}:=\ker \, t^{i+1}=\oplus_{j=0}^{i}\, Ke_j$.


\begin{lemma}\label{critN} 
Let $M$ be a $K[t]$-module. Suppose
 that $N(t,M)$ contains a Jordan module $N'$ such that $N'\supseteq \ker \, t$.
 Then $N(t,M)=N'$.
\end{lemma}

{\it Proof}. Since $N=\cup_{i\geq 0}  \, N_i$, it suffices to show that each
$N_i$ is contained in $N'$. We use induction on $i$. The case $i=0$ is true
 by the assumption. Suppose that $N_{i-1}\subseteq N'$ and let $u\in N_i$. Since
 $ tu\in N_{i-1}\subseteq N'$ and $tN'=N'$ ($N'$ is a Jordan module), we have
  $tu=tv$ for some $v\in N'$. Now, $t(u-v)=0$, hence $u-v=w\in N_0\subseteq N'$, and finally
$u=v+w\in N'$. It means that $N_i\subseteq N'$, as required. {\smallbox}

As we have seen in the Introduction, one can associate with a
non-scalar element $u$ of the Weyl algebra $A_1$ the inner
derivation $\ad \, u$ and the $\mathbb{N}$-filtered
 subspace $N(u)=N(u, A_1):= \cup_{i\geq 0}\, N(u,i,A_1)$ where
 $ N(u,i,A_1)=\ker (\ad \,u)^{i+1}$. The zero component of this filtration
is the centralizer $C(u,A_1)$ of the element $u$ in $A_1$. The vector space
 $N(u, A_1)$ is in fact a $\mathbb{N}$-graded algebra as follows from the formula,

\begin{equation}  
\label{dnab}
(\ad \, u)^n(ab)=\sum_{i=0}^n {n\choose i}\, (\ad \, u)^i(a)(\ad \, u)^{n-i}(b),
\end{equation}
where $a, b \in A_1$ and $n\geq 1$. The algebra $N(u,A_1)$ is a $K[\ad \, u]$-module,
 such that
$$
(\ad \, u)^i \, N(u, j, A_1)\subseteq N(u,j-i,A_1), \;\;
{\rm for \; all }\; i,j\geq 0, $$
 where we set $N(u, -i, A_1)=0$ for $i\geq 1$. Each $N(u, i, A_1)$ is a finitely generated
 $C(u, A_1)$-module (Proposition 10.2.(ii), \cite{Dix}). The associated graded algebra
of the $\mathbb{N}$-graded algebra $N(U,A_1)$,
$$ \CG (u,A_1):=\oplus_{i\geq 0}\, N(u,i,A_1)/N(u,i-1,A_1),$$
is a commutative domain (Proposition 10.2.(i), \cite {Dix} ).

{\bf Dixmier's Problem 4, \cite{Dix}}: {\em is the algebra} $\CG
(u, A_1)$
 {\em finitely generated?}

This problem is still open. In the next section a positive answer will be obtained
 for all homogeneous elements of the Weyl algebra. But first, we need a description of
the algebra $N(u,B)$ for an arbitrary homogeneous element of the algebra $B$.


\begin{theorem}\label{Nhomel} 
{\sc ( The Algebra $N(u,B)$ of a Homogeneous Element of the Algebra $B$)}
\begin{enumerate}
\item Let $u=\alpha X^n$ where $0\neq n \in \mathbb{Z}$ and  $\alpha $ is a nonzero
monic element of $K(H)$.

(i) The algebra $N(u,B)$ is generated by the algebra $C(u,B)$ and the element $H$.
 If $C(u,B)=K[v, v^{-1}]$ where $v=\beta X^t$ is chosen as in
Proposition \ref{CUB} then the algebra
 $N(u,B)$ is the skew Laurent extension $K[H][v, v^{-1}; \s^t]$. So, $N(u,B)$ is an
affine Noetherian algebra.

(ii) For $j\geq 0$, $N(u, j,B)=\sum_{i=0}^j \, C(u,B)H^i$ and
$H^j\in N(u,j,B)\backslash N(u,j-1,B)$.

(iii) The associated graded algebra $\CG (u,B)$ is the polynomial algebra
 $C(u,B)[h]$ with coefficients from $C(u,B)$ where $h:=H+C(u,B)\in N(u,1,B)\backslash
 C(u,B)$. Hence $\CG (u,B)$ is an affine commutative algebra.

\item Let $u\in K(H)\backslash K$. Then $N(u,B)=C(u,B)=K(H)$.
\end{enumerate}
\end{theorem}

{\it Proof}. 1. Clearly, $H\in N(u,1,B)\backslash C(u,B)$ since
$$ [u,H]=(\sn (H)-H)u=-nu$$
is a nonzero element of $C(u,B)$. Denote by $N'$ the subalgebra of $B$ generated by the algebra
$C(u,B)$ and the element $H$. The algebra $N'$ is a homogeneous subalgebra of the
 $\mathbb{Z}$-graded algebra $B=\oplus_{i\in \mathbb{Z}}\, K(H)X^i$ since $N'$ is
generated by the homogeneous elements $v$, $v^{-1}$ and $H$ of $B$. Using this fact we see that
$$ N'=K[H][v, v^{-1}; \s^t]$$
is a skew Laurent polynomial ring with coefficients from the polynomial ring
 $K[H]$. We aim to show that $N'=N(u,B)$. The inclusion
 $N'\subseteq N(u,B)$ is obvious since the algebra generators $v$, $v^{-1}$ and $H$
of $N'$ belong to $N(u,B)$. By Proposition \ref{CUB}.(1), $u=v^k$ where $k=t^{-1}n$ is
a natural number. The set of elements
\begin{equation}  
\label{vHnvj}
 \v_i(n^{-1}H)v^j, \; i\geq 0,\; j\in \mathbb{Z},
\end{equation}
is a $K$-basis of $N'$ where the polynomials $\v_i=\v_i(H)$ are defined in
(\ref{polynomialsvn}), and
$$\ad \, u\, (\v_i(n^{-1}H)v^j)=\v_{i-1}(n^{-1}H)v^{j+k}.$$
This means that $N'$ is a Jordan $K[\ad \,u]$-module which
contains $\ker \, \ad \, u=C(u,B)$ and $N'\subseteq N(u,B)$. By
Lemma \ref{critN}, $N'=N(u,B)$.

 It follows from (\ref{vHnvj}) that $N(u,j,B)=\sum_{i=0}^j\, C(u,B)\v_i(n^{-1}H)=
 \sum_{i=0}^j\, C(u,B)H^i$ and that $H^j\in N(u,j,B)\backslash N(u,j-1, B)$.
 So, we have proved the statements (i) and (ii). Statement (iii) is evident.

2. By the assumption, $u\in K(H)\backslash K$, thus, for each
nonzero $i \in \mathbb{Z}$, the element $\si (u)-u$ is nonzero
(since $K(H)^{\si }=K$, the algebra of $\si $-invariant elements
in the field $K(H)$). Let $w=\alpha X^m$ be a nonzero homogeneous
 element of $B$. For  $i\geq 1$,
$$(\ad \, u)^iw=(u-\s^m (u))^iw\neq 0.$$
It follows easily from this fact that $N(u,B)=K(H)$. {\smallbox}


\section{Centralizer and $N(u, A_1)$ of a Homogeneous Element of the Weyl Algebra}

In this section, for an arbitrary homogeneous element $u$ of the Weyl algebra $A_1$,
 algebra generators are found for
  the algebras $C(u,A_1)$ (Proposition
\ref{centralizerhomweyl}) and $N(u,A_1)$ (Theorem \ref{Nhomweyl}).
For certain homogeneous elements of $A_1$ their centralizers were
described in Proposition 5.3, \cite{Dix}. We shall see that the
Dixmier's Problem 4 has positive answer for all homogeneous
elements of the Weyl algebra
 (Theorem \ref{Nhomweyl}).

Consider an element $u=\alpha v_n\in A_1$ with $n\neq 0$ and a nonzero
monic polynomial $\alpha $
of $K[H]$. If $n>0$ then $u=\alpha X^n$, and if $n<0$ then $u=\alpha Y^n$. The element $u$ is
a monic
homogeneous element of the algebra $B$ since $\alpha $ is a monic polynomial and
\begin{equation}  
\label{Yn=nnXn}
Y^n=Y^nX^nX^{-n}=(-n,n)X^{-n}=H(H+1)\cdots (H+n-1)X^{-n}=(-1)^nn!\v_n X^{-n}.
\end{equation}
By Proposition \ref{CUB}.(1), $C(u, B)=K[v, v^{-1}]$ where $v=\beta X^t$ is the
canonical generator of the algebra $C(u, B)$, $0\neq \beta \in K(H)$ and the integer $t$ has the same sign
 as $n$. Moreover,
\begin{equation}  
\label{vm=u}
v^m=u, \;\; m=t^{-1}n\geq 1.
\end{equation}
The centralizer $C(u, A_1)=A_1\cap C(u,B)=\oplus_{i\in H}\, Kv^i$ where
 $H=\{i\in \mathbb{Z}: \, v^i\in A_1\}$. By \cite{Ami} or (Theorem 4.2, \cite{Dix}),
 the algebra $C(u,A_1)$ is a finitely generated $K[u]$-module. Since $u=v^m$ for some $m\geq 1$,
using the graded argument we have $H=\{i\geq 0: \, v^i\in A_1\}$.
For $i=0, 1, \ldots , m-1$, we denote by $\g_i$ a monic polynomial of $K[H]$ of minimal
possible degree, say $d_i$, such that $\g_iv^i\in A_1$. We denote by $\d $
the inner derivation $\ad \, u$ of the Weyl algebra $A_1$. Then
$$
A_1\ni \d^{d_i}(\g_iv^i)=(-n)^{d_i}d_i!v^iu^{d_i}=(-n)^{d_i}d_i!v^{i+d_im}.
$$
Thus, we can define the following non-negative integers,
$$
\mu_i:=\min \{ j\geq 0 \, | \, v^j\in A_1, \; j \equiv i\, ( {\rm mod} \, m) \},  \;\;
{\rm for \; \, each }\;\; i=0, 1, \ldots , m-1.
$$
Then
\begin{equation}  
\label{H=Umi} H=\cup_{i=0}^{m-1}\, \{ \mu_i +m \mathbb{N} \},
\end{equation}
a disjoint union. The next result describes the centralizer of an arbitrary
homogeneous element of the Weyl algebra $A_1$.

\begin{proposition}\label{centralizerhomweyl} 
 {\sc (Centralizer of a Homogeneous Element of the Weyl Algebra)}
\begin{enumerate}
\item Let $u=\alpha v_n\in A_1$ where $0\neq n\in \mathbb{Z}$ and
$\alpha $ is a monic polynomial of $K[H]$. Then
$C(u, A_1)=\oplus_{i=0}^{m-1}\, K[u]v^{\mu_i}$.
\item Let $u\in K[H]\backslash K$. Then $C(u, A_1)=K[H]$.
\end{enumerate}
\end{proposition}

{\it Proof}. $1.$  By (\ref{H=Umi}), $C(u, A_1)=\oplus_{j\in H}\,
Kv^j = \oplus \, \{ Kv^j\, |\, i=0, \ldots ,m-1 \; {\rm and} \;
j\in \mu_i +m\mathbb{N} \} =\oplus_{i=0}^{m-1}\, K[u]v^{\mu_i}$
since $v^m=u$.

$2.$ By Proposition \ref{CUB}.(2),  $C(u, B)=K(H)$, thus
$C(u, A_1)=A_1\cap C(u,B)=K[H]$. {\smallbox}

Observe that $X^t$ is  equal to $v_t$ if $t>0$, and to
$(t,-t)^{-1}v_t=(H(H+1)\cdots (H-t-1))^{-1}v_t$ if $t<0$ (by
(\ref{Yn=nnXn})). Thus the canonical generator $v$ can be written
in the form $\g v_t$ where
 $\g =\beta $ if $t>0$, and $\g =\beta (t,-t)^{-1}$ if $t<0$. The element $\g $ is
 a monic element of $K(H)$. Set $\mu :=\max \{ \mu_0, \ldots , \mu_{m-1}\}$.
Then
\begin{equation}  
\label{vjinA1}
v^i=\g \s^t(\g )\cdots \s^{(i-1)t}(\g ) v_{it}\in A_1, \;\; {\rm for \; all}\;
 i\geq \mu ,
\end{equation}
hence,
\begin{equation}  
\label{gginKH}
\g \s^t(\g )\cdots \s^{(i-1)t}(\g )\in K[H], \;\; {\rm for \; all}\;
 i\geq \mu .
\end{equation}
For each $i=1, \ldots , \mu -1$, there exists a unique monic polynomial $g_i\in K[H]$
 of minimal possible degree such
that $g_iv^i\in A_1$. The polynomial $g_i$ is the {\em denominator} of the rational
function $\g \s^t(\g )\cdots \s^{(i-1)t}(\g )$ multiplied by a proper nonzero scalar. By definition, the denominator of a rational function
$\alpha =pq^{-1}$ $(p,q\in K[H])$ is $q$ provided $\gcd (p,q)=1$.
 Clearly,
$$ K[H]v^i\cap K[H]v_{ti}=K[H]g_iv^i, \;\; i=1, \ldots , \mu -1.$$
It follows from the equality $v_{-k}v_k=(-k, k)\in K[H]$, $k\in
\mathbb{Z}$, that
$$ v_k^{-1}=(-k, k)^{-1}v_{-k}.$$
 Now, by (\ref{vjinA1}),
 $$v^{-i}=v^{-1}_{it}(\g \s^t(\g )\cdots \s^{(i-1)t}(\g ))^{-1}=
\{ (-it, it)\s^{-it}(\g \s^t(\g )\cdots \s^{(i-1)t}(\g ))\}^{-1}v_{-it}.$$
 By (\ref{gginKH}),
\begin{equation}  
\label{vit=vN}
v_{-it}=(-it, it)\s^{-it}(\g \s^t(\g )\cdots \s^{(i-1)t}(\g ))v^{-i}\in N(u, A_1), \;
 {\rm for \; all}\; i\geq \mu .
\end{equation}
For each $i\geq 1$, there exists a unique monic polynomial $f_i\in K[H]$ such that
$$ K[H]v^{-i}\cap K[H]v_{-it}=(f_i)v_{-it},$$
where $(f_i)=f_iK[H]$. By (\ref{vit=vN}), $f_i=1$ for all $i\geq \mu $. For each
 $i=1, \ldots , \mu-1 $, the polynomial $f_i$ is the denominator of the rational
function   $(-it, it)\s^{-it}(\g \s^t(\g )\cdots \s^{(i-1)t}(\g ))$ multiplied by a proper
nonzero scalar.

Let $R=\cup_{i\in \mathbb{N}}\, R_i$ be an $\mathbb{N}$-graded
algebra  and
 $\gr \, R=\oplus_{i\in \mathbb{N}}\, R_i/R_{i-1}$ be its associated graded algebra.
Denote by $ \pi :R\ra \gr \, R$ {\em the principal symbol map}
defined  $\pi (r)=r+R_{i-1}$ where $r\in R_i\backslash R_{i-1}$.

{\sc Definition}. A basis $E=\{ e_j, \, j \in J\}$ of the algebra $R$ is called a
{\em principal } basis iff the set $\pi (E)=\{ \pi (e_j), \, j\in J\}$ is a basis
of the associated graded algebra $\gr \, R$.

Suppose that $F=\{ f_l, \, l \in L\}$  is a basis of the algebra $\gr \, R$ such
 that each element $f_l$ is a homogeneous element of the algebra $\gr \, R$. For each $f_l$
 we fix its  preimage $e_l$ under the principal symbol map $\pi $, that is
$ \pi (e_l)=f_l$. Then the set $E=\{ e_l , \, l\in L\}$ is a principal basis of the
 algebra $R$.

The next theorem describes the algebra $N(u,A_1)$ for an arbitrary homogeneous element
 of the Weyl algebra $A_1$, and gives a positive answer to the Dixmier's Problem 4 for all such
 elements. In the next section, this result will lead us to a solution of the Dixmier's
Problem 5.


\begin{theorem}\label{Nhomweyl} 
{\sc (The Algebra $N(u,A_1)$ of a Homogeneous Element of the Weyl Algebra)}
\begin{enumerate}
\item Let $u=\alpha v_n \in A_1$ where $0\neq n \in \mathbb{Z}$ and $\alpha $ is a monic
polynomial of $K[H]$. We keep the notation above, then

(i) $N(u, A_1)=\oplus_{i\geq \mu}\, K[H]v_{-it}\, \oplus \,
\big( \oplus_{i=1}^{\mu -1}\,K[H]f_iv_{-it}\big) \oplus\, K[H]\oplus \,
 \big( \oplus_{i=1}^{\mu -1}\, K[H]g_iv^i \big) \oplus \,
\big( \oplus_{i\geq \mu}\, K[H]v^i\big) .$

(ii) The set $S=\{ H^kv_{-it}, \, H^kf_jv_{-jt}, \, H^k, \, H^kg_jv^j, \, H^kv^i
 \, | \, k\geq 0, i\geq \mu, j=1, \ldots , \mu -1 \}$  is a principal basis of the
 algebra $N(u, A_1)$, with
$\ndeg \, v_{-it}=i(|t|+ \deg_H\, \g)$ for $i\geq \mu $,
$\ndeg \, f_jv_{-jt}=\deg \, f_j+j(|t|+ \deg_H\, \g)$ and
$\ndeg \, g_jv^j=\deg_H\, g_j$ for $j=1, \ldots , \mu -1 $.

(iii) The algebra $\CG (u,A_1)$ is an affine (commutative) algebra, hence Noetherian.

(iv) The algebra $N(u,A_1)$ is an affine Noetherian algebra.

\item Let $u\in K[H]\backslash K$. Then $N(u,A_1)=K[H]=C(u,A_1)$.
\end{enumerate}
\end{theorem}

{\it Proof}. 1.(i).  By Theorem \ref{Nhomel}.(1), the algebra
$N(u,B)=\oplus_{i\in \mathbb{Z} } \, K[H]v^i $ and the Weyl
algebra $A_1=\oplus_{j\in \mathbb{Z} } \, K[H]v_j $  are
homogeneous subalgebras of the algebra $B=\oplus_{j\in \mathbb{Z}
} \, K[H]X^j $. So, the intersection
$$
N(u, A_1)= A_1\cap N(u,B)= \oplus_{i\in \mathbb{Z} } \,
(K[H]v^i\cap K[H]v_{it} ),$$ is a homogeneous subalgebra of the
Weyl algebra $A_1$. If we recall the definition of the polynomials
$f_i$ and $g_i$ then the result follows immediately from the fact
above and (\ref{vjinA1}), (\ref{vit=vN}).

(ii)  Since $N(u,A_1)$ is a homogeneous subalgebra of $A_1$, it is easy to see
 that the set $S$  is a principal basis of $N(u,A_1)$. Since
$\ndeg \, H=1= \deg_H \, H$ and
$$
v_{-it}=(-it, it)\s^{-it}(\g \s^t(\g )\cdots \s^{(i-1)t}(\g ))v^{-i}, \;
 {\rm for \; all}\; i\geq \mu , $$
we have $\ndeg \, v_{-it}=\deg_H\,(-it, it)\s^{-it}(\g \s^t(\g )\cdots
\s^{(i-1)t}(\g )) =i(|t|+ \deg_H\, \g)$. For each $j=1, \ldots , \mu -1 $,
$$
f_jv_{-jt}=f_j(-jt, jt)\s^{-jt}(\g \s^t(\g )\cdots \s^{(j-1)t}(\g ))v^{-j},
$$
hence,
$\ndeg \, f_jv_{-jt}=\deg_H\, f_j(-jt, jt)\s^{-jt}(\g \s^t(\g )\cdots
\s^{(j-1)t}(\g ))=\deg \, f_j+j(|t|+ \deg_H\, \g)$. The rest is obvious.

(iii)  Denote by $R$ the subalgebra of $\CG =\CG (u,A_1)$ generated by the principal
symbols of the elements $v_{-\mu t}$, $H$ and $v^\mu $. The set $S$ is a principal
 basis of the algebra $N(u,A_1)$, thus the set $\pi (S)=\{ \pi (s), \, s\in S\} $
is a basis of the algebra $\CG $. The algebra $\CG $ is affine since it is a finitely
 generated $R$-module  with  generators which are the principal symbols of the elements
$$ v_{-it}, \, f_jv_{-jt}, \, 1, \, g_jv^j,\; {\rm and} \; v^{it},
\;\; {\rm where }\;\;  i=\mu+1, \ldots , 2\mu -1; j=1, \ldots , \mu -1.$$

(iv) The algebra $N(u,A_1)$ is a Noetherian affine algebra since the
algebra $\CG $ is so.

2. By Theorem \ref{Nhomel}.(2), $N(u,A_1)=A_1\cap N(u, B)=A_1\cap K(H)=K[H]=
C(u,A_1)$. {\smallbox}

Let $u=\alpha v_n\in A_1$ be as in Theorem \ref{Nhomweyl}.(1). The algebra $N(u,A_1)$
 is a $\mathbb{Z}$-graded algebra with zero graded component $K[H]$, hence
 the (left and right) Krull dimension (in the sense of Rentschler and Gabriel, \cite{RG}),
\begin{equation}\label{KdimNKH}
\Kdim \, N(u,A_1)\geq \Kdim \, K[H]= 1.
\end{equation}

The homogeneous subalgebra $A$ of $N(u,A_1)$ generated by the homogeneous elements
 $y:=v_{-\mu t}$, $H$ and $x:=v^\mu $ is the generalized Weyl algebra
$$ A=K[H](\s^{\mu t} , a:=(-\mu t, \mu t) \s^{-\mu t}(\g \s^t(\g )\cdots \s^{(\mu -1)t}(\g ))),$$
since, by (\ref{vit=vN}),
$$ xH=\s^{\mu t}(H)x, \; yH=\s^{-\mu t}(H)y, \; yx=a, \; {\rm and} \; xy=\s^{\mu t} (a).
$$
In the terminology of \cite{Hod}, \cite{CB-H},  the algebra $A$ is called a {\em noncommutative
deformation of type-$A$ Kleinian singularity}. The algebra $A$ is a (left and right)
 Noetherian algebra, \cite{Bav}, \cite{Hod}.


\begin{corollary}\label{KdimN}  
Let $u=\alpha v_n\in A_1$ be as in Theorem \ref{Nhomweyl}.(1).
\begin{enumerate}
\item The algebra $N(u,A_1)$ is a finitely generated $A$-module.
\item The (left and right) Krull dimension of the algebra $N(u,A_1)$ is $1$.
\item If $\deg_H \, \alpha >0$ then the Weyl algebra $A_1$ is not a finitely
generated (left and right) $N(u,A_1)$-module.
\end{enumerate}
\end{corollary}

{\it Proof}. $1$ and $2$. By Theorem \ref{Nhomweyl}.(1), the algebra $N(u,A_1)$ is
a finitely generated (left and right) $A$-module, hence $\Kdim \, N(u,A_1) \leq
\Kdim \, A$. Since $a\neq 0$ and ${\rm char} \, K=0$, the  Krull dimension of
 the generalized Weyl algebra $A$ is $1$ (see \cite{Bav} or \cite{Hod}), hence
$\Kdim \, N(u,A_1)=1$, by (\ref{KdimNKH}).

Since $N(u, A_1)$ is a finitely generated $A$-module, the Weyl
algebra $A_1$
 is not a finitely generated $N(u, A_1)$-module iff it is not a
finitely generated $A$-module. So, it suffices to prove that $A_1$
is not
 a finitely generated $A$-module. By (\ref{vm=u}) and (\ref{vjinA1}), $\alpha v_n=u=
\g \s^t(\g ) \cdots \s^{(m-1)t}(\g )v_n$, hence $ 0<\deg_H \,
\alpha = \deg_H \, \g \s^t(\g ) \cdots \s^{(m-1)t}(\g )= m\deg_H
\, \g $, thus $\deg_H \, \g >0$ and $\deg_H \,a >0$. It suffices
to show that the factor module $A_1/A$ is not a finitely generated
$A$-module. Observe that $A$ is a homogeneous subalgebra
 of $A_1$, and that
$$M:=(\oplus_{i\in \mathbb{Z}}\, K[H]v_{i\mu t})/A=
\oplus_{i\geq 1}\; ( K[H]v_{i\mu t}/K[H]v^{i\mu })$$ is an
$\mathbb{N}$-graded $A$-submodule of $A_1/A$. The $i$'th component
of $M$, $M_i:=K[H]v_{i\mu t}/K[H]v^{i\mu }$, as a $K[H]$-module,
is canonically isomorphic to
$$K[H]v_{i\mu t}/K[H]\g_{i\mu }v_{i\mu t}\simeq K[H]/K[H]\g_{i\mu }, $$
 where $\g_{i\mu }:=\g \s^t(\g ) \cdots \s^{(i\mu -1)t}(\g )$. Since $\g_{i\mu }M_i=0$,
 for all $i\geq 1$, the $K[H]$-module $M_i$ is a $K[H]$-torsion module. Each finitely
generated $K[H]$-torsion module over the generalized Weyl algebra $A$ has finite length
 and Gelfand-Kirillov dimension $\leq 1$ (\cite{Bav, Bav2}).

Suppose that $M$ is a finitely generated $A$-module, then $\GK (M)= 1$  since
$\dim_K \, M=\infty $. The algebra $A$ is a {\em somewhat commutative } algebra, \cite{Hod}, hence
there exists a natural number $c$ such that
$$ \sum_{i=1}^n\, \dim \, M_i\leq cn \;\ {\rm for \; all}\;\; n>>0,$$
which contradicts  $\sum_{i=1}^n\, \dim \, M_i =\sum_{i=1}^n\,
\deg_H \, \g_{i\mu }= \sum_{i=1}^n\, i\mu \deg_H \, \g=\mu \deg_H
(\g ) \, \frac{n(n+1)}{ 2}$.
 Thus $M$ is not a finitely generated $A$-module. {\smallbox}


\section{ Solution to the Dixmier's Problem 5}

In this section we apply the results from the previous sections to show that the
 Dixmier's Problem 5 has {\em negative solution}.

\begin{lemma}\label{alphadegd} 
Let $\alpha $ be a monic polynomial of $K[H]$ of degree $d\geq 1$ and
$u=\alpha X\in A_1$. Then
\begin{enumerate}
\item $C(u,A_1)=K[u]$.
\item $N(u,A_1)=\oplus_{i\geq 1}\, K[H]Y^i \; \oplus \; \oplus_{i\geq 0}\, K[H]u^i $ and
the set $\{ \v_iY^{j+1}, \v_iu^j\, | \, i,j\geq 0\} $ is a basis of the algebra $N(u,A_1)$.
\item $Y\in N(u, d+1, A_1)\backslash N(u, d, A_1)$ and, for $k\geq 1$,
$N(u, k, A_1)= \oplus_{i,j\geq 0}\, \{ K\v_iY^j\, | \, i+(d+1)j\leq k\}\; \oplus \;
\oplus_{i=0}^k\, K[u]\v_i$.
\end{enumerate}
\end{lemma}

{\it Proof}. $1.$ By Proposition \ref{CUB}.(1), the algebra $C(u,B)=K[u,u^{-1}]$, hence, by
 Proposition \ref{centralizerhomweyl}.(1),  we have $C(u,A_1)=K[u]$.

$2.$ The element $u$ is the canonical generator of the algebra $C(u,B)=K[u,u^{-1}]$.
 Since $\alpha \in K[H]$, by Theorem \ref{Nhomweyl}.(1),
 $N(u,A_1)=\oplus_{i\geq 1}\, K[H]Y^i \; \oplus \; \oplus_{i\geq 0}\, K[H]u^i $. The rest is evident.

$3.$ By Theorem \ref{Nhomel}, we have $H\in N(u, 1 , A_1)\backslash N(u, 0 , A_1)$, hence
 $$Y=HX^{-1}=H\s^{-1}(\alpha )(\alpha X)^{-1}=H\s^{-1}(\alpha ) u^{-1}\in
N(u, d+1, A_1)\backslash N(u, d, A_1),$$
 and
$$\v_i (H)Y^j\in N(u, i+(d+1)j, A_1)\backslash N(u, i+(d+1)j-1, A_1), \;\;
{\rm  for \;\; all }\;\; i,j\geq 0.$$
 Now, the result follows from Statement 2. {\smallbox}

Let $u=\alpha X$ be as in Lemma \ref{alphadegd}. We denote by $\d $ the inner derivation
$\ad \, u$ of the Weyl algebra $A_1$.

For each $i\geq 1$, $\d (\v_i)=(\s (\v_i )-\v_i)u=\v_{i-1}u$, hence

\begin{equation}  
\label{divi=ui}
\d^i (\v_i)=u^i.
\end{equation}
Clearly,
\begin{eqnarray*}
Y^i& = &(HX^{-1})^i=H(H+1)\cdots (H+i-1)X^{-i}\\
&= & H(H+1)\cdots (H+i-1)\s^{-1}(\alpha )\s^{-2}(\alpha )\cdots
\s^{-i}(\alpha )u^{-i}\\
&=& (H^{i(d+1)}+\cdots )u^{-i}=(-1)^{i(d+1)}[i(d+1)]!\,
\v_{i(d+1)}u^{-i}+\cdots,
\end{eqnarray*}
 where by three dots we denote, as
usually, elements of smaller nilpotent degree. So,

\begin{equation}  
\label{dyi=uid}
\d^{(d+1)i}(Y^i)=(-1)^{i(d+1)}[(d+1)i]!\,u^{id}.
\end{equation}
Using (\ref{dnab}), we have

\begin{equation}  
\label{did1j} \d^{i+(d+1)j}(\v_iY^j)= {i+(d+1)j\choose i} \d^i
(\v_i)\d^{(d+1)j}(Y^j) =(-1)^{j(d+1)} {i+(d+1)j\choose
i}[(d+1)j]!\,u^{i+dj},
\end{equation}
for all $i,j\geq 0$.


\begin{corollary}\label{Dixpr} 
 {\sc (Solution to the Dixmier's Problem 5)}
 Let $u=\alpha X$ be as in Lemma \ref{alphadegd}. Then
 $I_k=u^{k-[\frac{k}{d+1}] }\, K[u]$,  for all $k\geq 1$. In particular,
$I_1=uK[u]$ and $I_{i(d+1)-1}=I_{i(d+1)}=u^{id}K[u]$, for all
$i\geq 1$. Hence, $I_1I_{i(d+1)-1}\neq I_{i(d+1)}$, for all $i\geq
1$, and the Dixmier's Problem 5
 has negative solution.
\end{corollary}

{\it Proof}. By Lemma \ref{alphadegd}.(3) and (\ref{did1j}),
 $I_k=u^{k-[\frac{k}{ d+1}] }\, K[u]$,  for all $k\geq 1$. The rest is obvious.
 {\smallbox}

It turns out that, for the element $u$ as above, the algebra $N(u,A_1)$ is a
generalized Weyl algebra of a special sort. So, applying the results of the papers
\cite{Bav}--\cite{Bav2}, \cite{Hod}, where these algebras were studied, we can say a lot
 about them. We collect some of the results in the following corollary.

\begin{corollary}\label{Nhomweylgwa}  
Let $u\in A_1$ be as in Lemma  \ref{alphadegd}. Then
\begin{enumerate}
\item The algebra $N(u,A_1)$ is a generalized Weyl algebra $K[H](\s , H\s^{-1}(\alpha
))$, a so-called  {\em noncommutative deformation of type
$A$-Kleinian singularity} in the terminology of \cite{Hod},
\cite{CB-H}.
\item The algebra $N(u,A_1)$ is simple iff, for any two distinct monic irreducible
factors $p$ and $q$ from $K[H]$ of the polynomial $H\s^{-1}(\alpha )$,
 there is no an integer $i$ such that $\si (p)=q$.
\item The algebra $N(u,A_1)$ has only finitely many (two-sided) ideals, they are
 classified in \cite{Bav1}. Each nonzero ideal has finite codimension in $N(u,A_1)$.
\item The Krull dimension of the algebra $N(u,A_1)$ is $1$.
\item Let $H\s^{-1}(\alpha )=p_1^{n_1}\cdots p_s^{n_s}$ be a product of distinct
monic irreducible polynomials. The global dimension

\[ {\rm gl.dim}\, N(u,A_1) =\left\{ \begin{array}{llll}
\infty & \mbox{, if there exists $n_i\geq 2$; }\\
2 & \mbox{, if $n_1=\cdots =n_s=1$ and $\si (p_j)=p_k$    }\\
{} & \mbox{\,\, for some $j\neq k$ and some integer $i$; }\\
1 & \mbox{, otherwise.}
\end{array}
\right.
\]

\end{enumerate}
\end{corollary}

{\it Proof}. $1.$ By Lemma \ref{alphadegd}.(2), the algebra
$N(u,A_1)$ is generated by the elements $Y$, $H$ and $X'=\alpha
X$. Since
$$
X'H=\s (H)X', \; YH=\s^{-1}(H)Y,\; YX'=H\s^{-1}(\alpha )\; {\rm and}\; X'Y=\s (H\s^{-1}(\alpha )),
$$
the algebra $N(u,A_1)$ is isomorphic to the generalized Weyl algebra
$K[H](\s , H\s^{-1}(\alpha ))$ in a view of the decomposition from Lemma
 \ref{alphadegd}.(2).

$2$ and $3$. These results were proved in \cite{Bav, Bav1, Bav2}.

$4$ and $5$. These results were proved in \cite{Bav, Bav2, Hod}. {\smallbox}


\begin{corollary}\label{weylN}  
Let $u\in A_1$ be as in Lemma  \ref{alphadegd}.
 Then the left $N(u,A_1)$-module $M=A_1/N(u,A_1)$ is a $K[H]$-torsion,
 not finitely generated left  $N(u,A_1)$-module of Gelfand-Kirillov dimension $1$.
Each finitely generated $N(u,A_1)$-submodule of $M$ has finite length. The set of
isomorphism classes of simple subfactors of all finitely generated
 $N(u,A_1)$-submodules of $M$is a finite set.
\end{corollary}

{\it Proof}.
The algebra $N(u,A_1)$ is a homogeneous subalgebra of the Weyl
algebra $A_1$, thus the $N(u,A_1)$-module $M=\oplus_{i\geq 1}\,
M_i$ is an $\mathbb{N}$-graded $N(u,A_1)$-module where the $i$'th
component $M_i$, as a $K[H]$-module, is canonically isomorphic to
$$K[H]/u^iK[H]= K[H]/\alpha_iX^iK[H]\simeq K[H]/\alpha_iK[H],$$
where $\alpha_i:=\alpha \s(\alpha )\cdots \s^{i-1} (\alpha )$.
Since $\alpha_i M_i=0$, for all $i\geq 1$, the module $M$ is a
$K[H]$-torsion module. Each finitely generated $K[H]$-torsion
module over a generalized Weyl algebra of the type $K[H](\s ,
a\neq  0)$, for example $N(u,A_1)$, has finite length and
Gelfand-Kirillov dimension $\leq 1$, thus $\GK (M)\leq 1$. Observe
that the $N(u,A_1)$-submodule $L$ of $M$ generated by
 the element $\bar{X}=X+N(u,A_1)$ is not finite dimensional since
$$ u^i\bar{X}=\alpha_iX^i\bar{X}=X\s^{-1}(\alpha_i)X^i +N(u,A_1)=
X(\s^{-1}(\alpha_i)-\alpha_i)X^i+N(u,A_1)=(\alpha_i-\s (\alpha_i))\bar{X^i}\neq 0$$
since $0<\deg (\alpha _i -\s (\alpha_i))<\deg \, \alpha_i$. So, $1\leq \GK (L)\leq \GK (M)\leq 1$, hence
 $\GK (M)=1$.

A finitely generated $N(u,A_1)$-submodule, say $V$, of $M$ is a submodule of the
module $U_s$ generated by $\oplus_{i=1}^s\, M_i$ for some $s$. The $N(u,A_1)$-module
 $U_s$ is an epimorphic image of the $N(u,A_1)$-module
$\oplus_{i=1}^s\, N(u,A_1)/N(u,A_1)\alpha_i$.
 Each $N(u,A_1)$-module $N(u,A_1)/N(u,A_1)\alpha_i$ has finite length, and the set of all
isomorphic classes of all simple subfactors of all modules $N/N\alpha_i$ is a finite
set since $\alpha_i:=\alpha \s(\alpha )\cdots \s^{i-1} (\alpha )$ (see \cite{Bav, Bav2} for details).
 Now the result follows. {\smallbox}


\section{Classification of Homogeneous Elements of the Weyl Algebra and
 the Dixmier's Problem 4}

Let $u$ be a non-scalar element of the Weyl algebra $A_1$. The corresponding
inner derivation $\ad \, u$ of $A_1$ is denoted by $\d $. Denote by
$\Ev (\ad \, u, A_1)=\Ev (u,A_1)$ the set of all eigenvalues of the linear
 map $\d $ acting in the vector space $A_1$. For an eigenvalue $\l $
 of $\d $ we denote by $D(u, \l ,A_1)$ the set of all eigenvectors of $\d $
 with the eigenvalue $\l $. The map $\d $ is a derivation of the Weyl algebra $A_1$,
 so the set $\Ev (u, A_1)$ is an additive submonoid
 of the field $K$, and the vector space
$$ D(u)^{ev}=D(u, A_1)^{ev}:=\oplus_{\l \in \Ev(u,A_1)}\, D(u,\l , A_1)$$
is a $\Ev (u,A_1)$-graded algebra, that is,
$$D(u,\l , A_1)D(u,\mu , A_1)\subseteq D(u,\l +\mu, A_1), \;\; {\rm for \;\; all}\;\;
 \l ,\mu \in \Ev (u,A_1).$$
Let a field $\bK $ be an algebraic closure of the field $K$. The tensor product of algebras
 $\bA_1=K\t A_1$ over the field $K$ is the Weyl algebra over the field $\bK $ which
 contains the Weyl $K$-algebra $A_1$. Then
$$ D(u)=D(u,A_1):=A_1\cap D(u, \bA_1)^{ev}$$
is a $K$-subalgebra of $A_1$ which contains the algebra
$D(u)^{ev}$ but does not necessarily coincide with this algebra.
The next result, obtained by Dixmier, \cite{Dix}, classifies
non-scalar elements of the Weyl algebra $A_1$ with respect to the
properties of the corresponding inner derivations of elements.

\begin{theorem}\label{Dixclasel} 
{\sc (The Dixmier's Classification of Non-scalar Elements of the
Weyl Algebra)}
 The set of non-scalar elements of the Weyl algebra $A_1$ is a disjoint
union of the following subsets:

(i) $\D_1 =\{x\in A_1\backslash K:\,  N(x)=A_1, \;\; D(x)=C(x)\}.$

(ii) $\D_2 =\{x\in A_1\backslash K:\,  N(x)\neq A_1, \;\; N(x)\neq C(x), \;\; D(x)=C(x)\}.$

(iii) $\D_3 =\{x\in A_1\backslash K:\,  D(x)=A_1, \;\; N(x)=C(x)\}.$

(iv) $\D_4 =\{x\in A_1\backslash K:\,  D(x)\neq A_1, \;\; D(x)\neq C(x), \;\; N(x)=C(x)\}.$

(v) $\D_ 5=\{x\in A_1\backslash K:\,  D(x)=N(x)=C(x) \}.$ {\smallbox}
\end{theorem}
Each subset $\D_i$ of $A_1$ is a non-empty set.

{\sc Definition}. Elements of $\D_3\cup\D_4$ (resp. of $\D_3$) are called
 elements of {\em semi-simple type} (resp. of {\em strongly semi-simple type}).
 Dixmier classified elements of strongly semi-simple type, Theorem 9.2,
\cite{Dix}: $x\in \D_3$ iff there exists an automorphism $\tau \in
\Aut_K \, A_1$ such that $\tau (x)=\l Y^2+\mu X^2+\nu $ for some
scalars
 $\l , \mu $ and $\nu $ such that $\l \neq 0$ and $\mu \neq 0$. It can be easily seen
 that, if the polynomial $\l t^2+\mu $ has a root in the field $K$ then there exists
 an automorphism $\tau_1 \in \Aut_K \, A_1$ such that $\tau_1 (x)=\alpha H+\beta $
for some scalars $0\neq \alpha$ and $\beta $ (see Corollary 9.3, \cite{Dix}).


\begin{theorem}\label{Clashomel} 
{\sc (Classification of Homogeneous Elements of the Weyl Algebra)}
 Let $u=\alpha v_i$, $\alpha \in K[H]$, $i\in \mathbb{Z}$, be a homogeneous non-scalar
 element of $A_1$. Then $u\in {\Delta}_1\cup {\Delta}_2\cup {\Delta}_3\cup
{\Delta}_5$. In more detail,
\begin{enumerate}
\item $u\in {\Delta}_1 $ $\Leftrightarrow $ $\alpha \in K^*$ and $i\neq 0$.
\item $u\in {\Delta}_2 $ $\Leftrightarrow $ $\alpha \not\in K$ and $i\neq 0$.
\item $u\in {\Delta}_3 $ $\Leftrightarrow $ $\deg_H\, \alpha =1$ and $i=0$.
\item $u\in {\Delta}_5 $ $\Leftrightarrow $ $\deg_H\, \alpha >1$ and $i=0$.
\end{enumerate}
\end{theorem}

{\it Proof}. If $\alpha \in K^*$ and $i\neq 0$ then by \cite{Dix}, Proposition 10.3,
$N(u,A_1)=N(v^{|i|}_{\pm 1}, A_1)=N(v_{\pm 1}, A_1)=A_1$, hence $u\in {\Delta}_1 $.
 If $\alpha \not\in K$ and $i\neq 0$ then $u\in {\Delta}_2 $ , by
Theorem \ref{Nhomweyl}.(1). If $\deg_H\, \alpha =1$ and $i=0$, i.e. $u=\l H+\mu$
for some scalars $\l \neq 0$ and $\mu $, then $u\in {\Delta}_3 $. If
 $\deg_H\, \alpha >1$ and $i=0$ then, by Theorem \ref{Nhomweyl}.(2),
 $N(u,A_1)=C(u,A_1)=K[H]$. The element $u=\alpha $ is a homogeneous element of the
algebra $A_1$, thus the algebra $D(u, A_1)$ is a homogeneous subalgebra of $A_1$.
Suppose that $D(u, A_1)\neq C(u,A_1)$ then there exists a homogeneous element $\beta v_m$
$(\beta \in K[H])$ of $A_1$ and a nonzero scalar $\l $ such that
$\beta v_m\in D(u,\l , A_1)$, hence
 $\l \beta v_m=[\alpha , \beta v_m]=(\alpha -\s^m (\alpha ))\beta v_m$.
 So, $\l =\alpha -\s^m (\alpha )$, hence $\deg_H\, \alpha =1$, a contradiction. This means
that $D(u, A_1)=C(u, A_1)$ and $u\in {\Delta}_5 $. This finishes the proof of the theorem.
 {\smallbox}

\begin{corollary}\label{Dixpr4hom}
Let $u$ be a homogeneous element of weakly nilpotent type of the Weyl
 algebra $A_1$, i.e. $u\in \D_2$.
 Then
\begin{enumerate}
\item The associated graded algebra $\CG (u,A_1)$ is an affine commutative algebra
 (thus, the Dixmier's Problem 4 has a positive answer for homogeneous elements of nilpotent type)
 hence the algebra $N(u,A_1)$ is an affine Noetherian algebra.
\item The algebra $A_1$ is not a finitely generated (left and right) $N(u,A_1)$-module.
\end{enumerate}
\end{corollary}

{\it Proof}. $1$. By Theorem \ref{Clashomel}.(2), each homogeneous element of nilpotent type
of the Weyl  algebra $A_1$ has the form as in Theorem \ref{Nhomel}.(1). Now observe that
the result is already was proved in Theorem \ref{Clashomel}.(2).(iii) and (iv).

$2$. This follows from Corollary \ref{KdimN}.(3) and Theorem \ref{Clashomel}.(2).
 {\smallbox}

Let $a$ and $p$ be nonzero elements of the Weyl algebra $A_1$
satisfying $[a,p]=\l p$ for some $0\neq \l \in K$. Then the
element $p$ is of nilpotent type, and
 $[a,cp]=\l cp$ for all nonzero
 elements $c\in C(a, A_1)$. The next result shows how the type of the element $p$
 changes when $p$ is multiplied by an element from $C(a,A_1)$.


\begin{corollary}\label{p=eigenval=cp}
Let $a\in \D_3(A_1)$ and $[a,p]=\l p$ for some $0\neq \l \in K$ and $p\in A_1$. Then
 $C(a,A_1)=K[a]$.
\begin{enumerate}
\item Suppose that $p\in \D_1(A_1)$ and $0\neq \alpha (t)\in K[t]$. Then

(i) $\alpha (a)p\in \D_1(A_1)$ if and only if $\alpha \in K^*$.

(ii) $\alpha (a)p\in \D_2(A_1)$ if and only if $\alpha \not\in K$.

\item Suppose that $p\in \D_2(A_1)$. Then $\alpha (a)p\in \D_2(A_1)$ for
 all nonzero polynomials $\alpha (t)\in K[t]$.
\end{enumerate}
\end{corollary}

{\it Proof}. The fact that $C(a,A_1)=K[a]$ easily follows from Theorem 9.2 and Corollary
9.4, \cite{Dix}.

We may assume that $K$ is an algebraically closed field. Then, by  Theorem 9.2 and Corollary
9.4, \cite{Dix}, there exists an automorphism $\nu \in \Aut_K(A_1)$ such that
$\nu (a)=\mu H+\g $ for some scalars $\mu \neq 0$ and $\g $. So, multiplying the element
 $a$ by $\mu^{-1}$ and adding an appropriate scalar to $H$, without loss of generality we may assume
 that $a=H$. Now, $p$ is a homogeneous nonscalar element of the algebra $A_1$, and the result follows from
 Theorem \ref{Clashomel}. {\smallbox}

\begin{center}
{\sc Acknowledgment }
\end{center}
\noindent The author would like to thank J. Dixmier and T. Lenagan
for comments.

Department of Pure Mathematics

University of Sheffield

Hicks Building

Sheffield S3 7RH

email: v.bavula@sheffield.ac.uk

\end{document}